\documentclass{amsart}

\theoremstyle{plain}
\newtheorem{theorem}{Theorem}[section]
\newtheorem{lemma}[theorem]{Lemma}
\newtheorem{proposition}[theorem]{Proposition}
\newtheorem{corollary}[theorem]{Corollary}

\theoremstyle{definition}
\newtheorem{defn}[theorem]{Definition}
\newtheorem{example}[theorem]{Example}

\theoremstyle{remark}
\newtheorem{remark}[theorem]{Remark}

\numberwithin{equation}{section}

\newcommand{\Aa}{{\mathbb A}}

\newcommand{\Pp}{{\mathbb P}}
\newcommand{\Zz}{{\mathbb Z}}
\newcommand{\Rr}{{\mathbb R}}

\newcommand{\SV}{{\mathcal S}}               
\newcommand{\PV}{{\mathcal V}}               
\newcommand{\PS}{{\mathcal X}}               
\newcommand{\VO}{{\mathcal V}^{\circ}}       

\newcommand{\SVA}{\tilde{\mathcal S}}        
\newcommand{\PVA}{\tilde{\mathcal V}}
\newcommand{\PSA}{\tilde{\mathcal X}}
\newcommand{\VAO}{\tilde{\mathcal V}^{\circ}}

\newcommand{\D}{{\mathcal D}}                
\newcommand{\I}{{\mathcal I}}                
\newcommand{\A}{{\mathcal A}}                
\newcommand{\B}{{\mathcal B}}

\newcommand{\LRM}[1]{{\mathcal M}(#1)}       
\newcommand{\SRM}[1]{{\mathcal M}^s(#1)}     

\newcommand{\refineq}{\preceq}
\newcommand{\refinedbyeq}{\succeq}

\newcommand{\codim}{\mathop{{\rm codim}}\nolimits}
\newcommand{\Coup}{\mathop{\rm Cpl}\nolimits}
\newcommand{\Gr}{\mathop{\rm Gr}\nolimits}
\newcommand{\Spec}{\mathop{\rm Spec}\nolimits}
\newcommand{\supp}{\mathop{\rm supp}\nolimits}
\newcommand{\val}{\mathop{\rm val}\nolimits}

\newcommand{\Case}[2]    {\noindent {\em Case~#1:~#2\/}}
\newcommand{\de}[2]{[#1,#2]}                          
\newcommand{\defterm}[1]{{\it #1\/}}
\newcommand{\fld}{{\bf k}}
\newcommand{\isom}{\cong}
\newcommand{\mm}{{\bf m}}
\newcommand{\Pic}{{\bf P}}
\newcommand{\sgn}{\varepsilon}
\newcommand{\sm}{\setminus}
\newcommand{\st}{~\mid~}       
\newcommand{\TP}{{\tau}}       
\newcommand{\ue}[2]{\{#1,#2\}} 
\newcommand{\x}{\times}

\begin{document}

\title{Geometry of graph varieties}
\author{Jeremy L. Martin}
\address{Department of Mathematics,
University of California, San Diego,
La Jolla, CA 92093-0112}
\curraddr{School of Mathematics,
University of Minnesota,
Minneapolis, MN 55455}
\email{martin@math.umn.edu}

\begin{abstract}
A {\em picture} $\Pic$ of a graph $G = (V,E)$ consists of a point $\Pic(v)$ for each vertex $v
\in V$ and a line $\Pic(e)$ for each edge $e \in E$, all lying in the projective plane over a
field $\fld$ and subject to containment conditions corresponding to incidence in $G$.  A {\em
graph variety} is an algebraic set whose points parametrize pictures of $G$.  We consider
three kinds of graph varieties: the {\em picture space} $\PS(G)$ of all pictures;  the {\em
picture variety} $\PV(G)$, an irreducible component of $\PS(G)$ of dimension $2|V|$, defined as
the closure of the set of pictures on which all the $\Pic(v)$ are distinct; and the {\em slope
variety} $\SV(G)$, obtained by forgetting all data except the slopes of the lines $\Pic(e)$.  
We use combinatorial techniques (in particular, the theory
of {\em combinatorial rigidity}) to obtain the following
geometric and algebraic information on these varieties:
	\begin{enumerate}
	\item A description and combinatorial interpretation of equations defining each 
		variety set-theoretically.
	\item A description of the irreducible components of $\PS(G)$.
	\item A proof that $\PV(G)$ and $\SV(G)$ are Cohen-Macaulay when $G$ satisfies a
		sparsity condition, {\it rigidity independence}.
	\end{enumerate}
In addition, our techniques yield a new proof of the equality of two matroids studied in 
rigidity theory.
\end{abstract}

\maketitle

\section{Introduction} \label{intro-section}

This paper initiates the study of certain algebraic varieties that parametrize plane pictures
$\Pic$ of a given graph $G$, with vertices $v$ and edges $e$ represented respectively by
points $\Pic(v) \in \Pp^{2}$ and lines $\Pic(e)$ connecting them in pairs.  Three such
varieties naturally arise.  First of all, there is the {\it picture space} $\PS(G)$ of all
pictures of $G$.  Usually, $\PS(G)$ is not irreducible.  It is therefore natural to restrict
attention to a second variety, namely the irreducible component of $\PS(G)$ containing as a
dense set those pictures in which the points $\Pic(v)$ are all distinct.  This most generic
component of the picture space is called the {\it picture variety} $\PV(G)$.  As we shall see
in Theorems~\ref{tree-poly-thm} and \ref{defining-eqns-thm}, $\PV(G)$ is cut out in $\PS(G)$
purely by equations relating the slopes of the lines $\Pic(e)$. The crucial matter for the
whole study is to understand the relations among these slopes.  This leads us to consider the
{\it slope variety} $\SV(G)$, which is essentially the projection of $\PV(G)$ on coordinates
$m_{e}$ giving the slopes of the lines $\Pic(e)$.

In a sequel to this paper~\cite{JLM2}, we study intensively the case where $G$ is the complete
graph $K_{n}$.  There we obtain very precise results, including the proof for $K_{n}$ of some
conjectures mentioned below, along with remarkable connections to the combinatorics of
matchings and planar trees.  Note that the problem of describing the slope variety
$\SV(K_{n})$ is of a very classical kind: it is exactly the problem of determining all
relations among the slopes of the $\binom{n}{2}$ lines connecting $n$ general points in the
plane.

Here we consider the features of varieties associated with an arbitrary graph $G=(V,E)$.  In
Section~\ref{general-section}, we define precisely the varieties $\PS(G)$, $\PV(G)$ and
$\SV(G)$ and give a fundamental construction, the decomposition of the picture space into
locally closed pieces which we call \defterm{cellules}.  In Section~\ref{rigidity-section}, we
show that the {\it generic rigidity matroid} $\LRM{V}$ studied by Laman {\it et.~al.}
(see \cite{Lam70},~\cite{GSS93}) appears as the algebraic dependence matroid of the slopes
(Theorem~\ref{matroid-thm}).

The heart of the paper is Section~\ref{main-section}.  For each set of edges forming a circuit
in the matroid $\LRM{V}$, we can write down an explicit determinantal formula for the
essentially unique polynomial relation among the corresponding slopes $m_{e}$.  The slope
relation induced by each circuit in $\LRM{V}$ turns out to be a remarkable polynomial.  All
its terms are squarefree, and they have a surprising combinatorial interpretation in terms of
decompositions of the given circuit into pairs of complementary spanning trees
(Theorem~\ref{tree-poly-thm}).  Moreover, we prove that precisely these relations cut out
$\PV(G)$ in $\PS(G)$ set-theoretically (Theorem~\ref{defining-eqns-thm}).

In Section~\ref{geom-section}, we show how the full component structure of $\PS(G)$ can be
economically described in terms of the rigidity matroid (Theorem~\ref{cellule-closure-thm}),
and show that when $\PS(G)=\PV(G)$, this variety has Cohen-Macaulay singularities
(Theorem~\ref{CM-thm}). We conjecture that the slope relations should cut out $\PV(G)$
scheme-theoretically as well as set-theoretically. We further suspect that they may always
form a {\it universal Gr\"obner basis} for the ideal of the slope variety, and moreover, that
both $\SV(G)$ and $\PV(G)$ are always Cohen-Macaulay.

When we first embarked upon the study of graph varieties, before obtaining the results
indicated above, we already had some reasons to think they might be interesting.  Since these
reasons remain relevant, let us mention them briefly.  Graph varieties are a particular class
of {\it configuration varieties}---subvarieties in a product of Grassmannians defined by
containment conditions among various subspaces of a fixed space.  Other examples of
configuration varieties are Bott-Samelson-Demazure varieties and somewhat more general
varieties introduced and studied by Magyar \cite{Mag97}.  The latter have very special
geometric properties, and it is natural to inquire to what extent these are shared by more
general configuration varieties. Graph varieties provide the simplest non-trivial examples not
fitting into Magyar's framework.  Furthermore, for $G=K_{n}$, the picture variety $\PV(G)$ is
a blowdown of the Fulton-Macpherson {\it compactification of configuration space} \cite{FM94},
which desingularizes it.  For general $G$, the same relation holds between $\PV(G)$ and the
DeConcini-Procesi {\it wonderful model of subspace arrangements} \cite{DP95}.  We expect that
$\PV(G)$ should not only be Cohen-Macaulay but should have rational singularities.  This would
be equivalent to a cohomology vanishing theorem for certain line bundles on the wonderful
model, raising an important question for further study.

This paper is essentially Chapter 1 of the author's Ph.D. dissertation~\cite{JLMPHD}, written
under the supervision of Mark Haiman.  The author wishes to thank Prof. Haiman for his advice
and encouragement.

\section{Definitions} \label{definition-section}

We work over an algebraically closed field $\fld$.  Affine and projective $n$-space over
$\fld$ are denoted by $\Aa^n$ and $\Pp^n$ respectively.

A \defterm{graph} is a pair $G=(V,E)$, where $V=V(G)$ is a finite set of \defterm{vertices}
and $E=E(G)$ is a set of \defterm{edges}, or unordered pairs of distinct vertices $\ue{v}{w}$.  
We frequently abbreviate $\ue{v}{w}$ by $vw$ when no confusion can arise (for instance, when
the vertices are one-digit positive integers). The vertices $v,w$ are called the
\defterm{endpoints} of the edge $vw$.  A \defterm{subgraph} of $G$ is a graph $G'=(V',E')$
with $V' \subset V$ and $E' \subset E$.  We define
	\begin{equation*}
	K(V') = \{vw \st v,w \in V', ~ v \neq w\}
	\end{equation*}
and
	\begin{equation*}
	E(V') =  E \cap K(V').
	\end{equation*}
The \defterm{complete graph on} $V$ is the graph $(V,K(V))$.  We write $K_n$ for 
the complete graph on $\{1, \dots, n\}$.  For $E' \subset E$ and $v \in V$, the
\defterm{valence} of $v$ with respect to $E'$ is
	\begin{equation*}
	\val_{E'}(v) = \big\vert \left\{ e \in {E'} \st v \in e \right\} \big\vert
	\end{equation*}
and the \defterm{vertex support} of $E'$ is
	\begin{equation*}
	V(E') = \left\{ v \in V \st \val_{E'}(v) > 0 \right\}.
	\end{equation*}
For $v_1, \dots, v_s \in V$, we set
	\begin{equation*}
	(v_1, \dots, v_s) = \left\{ v_1v_2,~ v_2v_3,~ \dots,~ v_{s-1}v_s \right\} 
	\subset E.
	\end{equation*}
If the $v_i$ are all distinct, then $(v_1, \dots, v_s)$ is called a
\defterm{path}.  If $v_1, \dots, v_{s-1}$ are distinct and $v_1=v_s$, then $(v_1,
\dots, v_s)$ is called a \defterm{polygon} or \defterm{$(s-1)$-gon}.  A polygon is more
commonly called a ``cycle'' or ``circuit,'' but we wish to reserve these words for other
uses.

A graph $G=(V,E)$ is \defterm{connected} if every pair of vertices are joined by a path, and
is a \defterm{forest} if at most one such path exists for every pair.  A connected forest is
called a \defterm{tree}.  A \defterm{spanning tree} of $G$ (or of $V$) is a tree $T \subset E$
with $V(T)=V$.  A \defterm{connected component} of $G$ is a maximal connected subgraph; every
graph has a unique decomposition into connected components (where some components may be
isolated vertices).

A \defterm{partition} of a finite set $V$ is a set $\A = \{A_1, \dots, A_s\}$ of pairwise
disjoint subsets of $V$ whose union is $V$.  The sets $A_i$ are called the \defterm{blocks} of
$\A$.  We write $\sim_\A$ for the equivalence relation on $V$ whose equivalence classes are
the blocks of $\A$.  We distinguish two extreme cases: the \defterm{discrete partition}
$\D_V$, all of whose blocks are singletons, and the \defterm{indiscrete partition} $\I_V$,
which has only one block.  Finally, if $\A$ and $\B$ are partitions of $V$, then we say that
$\A$ \defterm{refines} $\B$, written $\A \refineq \B$, if every block of $\A$ is contained in
some block of $\B$.  It is elementary that refinement is a partial ordering.

\section{The Picture Space and Picture Variety of a Graph} \label{general-section}

Throughout this section, we consider a graph $G = (V,E)$ with $|V| = n$ and $|E| = r$.  
Define
	\begin{equation} \label{grg}
	\Gr(G) = \prod_{v \in V} \Pp^2 ~\x~ \prod_{e \in E} \hat{\Pp}^2,
	\end{equation}
where $\hat{\Pp}^2$ denotes the dual projective plane of lines in $\Pp^2$.

For $\Pic \in \Gr(G)$, $v \in V$, and $e \in E$, we write $\Pic(v)$ and $\Pic(e)$ for
the projections of $\Pic$ on the indicated factors in (\ref{grg}).

\begin{defn}
A \defterm{picture} of $G$ is a point $\Pic \in \Gr(G)$ such that
$\Pic(v) \in \Pic(e)$ whenever the vertex $v$ is an endpoint of the edge $e$.
The \defterm{picture space} $\PS(G)$ is the set of all pictures of $G$.
\end{defn}

Note that $\PS(G)$ is Zariski closed in $\Gr(G)$, since the incidence conditions
may be expressed in terms of the Pl\"ucker coordinates.  Note also that if $G_1, \dots,
G_s$ are the connected components of $G$, then
	\begin{equation*}
	\PS(G) \isom \PS(G_1) \x \dots \x \PS(G_s).
	\end{equation*}

The equations defining $\PS(G)$ in homogeneous coordinates are awkward to work with
explicitly.  However, all the geometric information we will require about $\PS(G)$ can be
recovered from the following affine open subset of it, on which the defining equations
assume a more manageable form.

\begin{defn}
Fix homogeneous coordinates $[a_0 : a_1 : a_2]$ on $\Pp^2$, identifying $\Aa^2$ with the points
for which $a_0 \neq 0$ and giving $x = a_1/a_0$, $y = a_2/a_0$ as affine coordinates on $\Aa^2$.
The \defterm{affine picture space} $\PSA(G)$ is the open subvariety of $\PS(G)$ consisting of
pictures $\Pic$ such that all points $\Pic(v)$ lie in $\Aa^2$ and no line $\Pic(e)$ is 
parallel to the $y$-axis.
\end{defn}

Note that $\PSA(G)$ is open and dense in $\PS(G)$, and that $\PS(G)$ is covered by finitely
many copies of $\PSA(G)$.  In addition $\PSA(G)$ has affine coordinates
	\begin{equation*}
	\{x_v,y_v \st v \in V\} ~\cup~ \{m_e,b_e \st e \in E\},
	\end{equation*}
where $m_e$ and $b_e$ denote respectively the slope and $y$-intercept of the line
$\Pic(e)$.  Thus $\PSA(G)$ is the vanishing locus (in $\Aa^{2n+2r}$, identified with an open
subset of $\Gr(G)$) of the ideal generated by the $2r$ equations
	\begin{equation} \label{affineone}
	\begin{array}{c}
	y_v = m_e x_v + b_e, \\
	y_w = m_e x_w + b_e,
	\end{array}
	\end{equation}
for each edge $e = vw$.  Eliminating the variables $b_e$ from ({\ref{affineone})
produces $r$ equations
	\begin{equation} \label{affinetwo}
	(y_v - y_w) = m_e(x_v-x_w).
	\end{equation}
We may also eliminate the variables $y_v$.  For each polygon $P = (v_1,
\dots, v_s, v_1)$ of $G$, we sum the equations (\ref{affinetwo}) over the edges of $P$,
obtaining the equation
	\begin{equation} \label{affinethree}
	L(P) = \sum_{i=1}^s m_{e_i} (x_{v_i} - x_{v_{i+1}}) = 0.
	\end{equation}
where $e_i=v_iv_{i+1}$ and the indices are taken modulo $s$. Given a
solution in $\{m_e,x_v\}$ of the equations (\ref{affinethree}), we may choose one
$y$-coordinate arbitrarily and use (\ref{affineone}) and (\ref{affinetwo}) to recover the
coordinates $y_v$ and $b_e$.  Thus, putting
	\begin{equation*}
	R'_G = \fld[m_e,x_v \st e \in E, ~ v \in V],
	\end{equation*}
we see that $\PSA(G) \cong \Aa^1 \x X$, where $X$ is the subscheme of $\Spec R'_G 
\isom \Aa^{|V|+|E|}$ defined set-theoretically by the equations (\ref{affinethree}).
We will show eventually that in fact, the equations defining the varieties $\PVA(G)$
and $\SVA(G)$ lie in the subring
	\begin{equation*}
	R_G = \fld[m_e \st e \in E].
	\end{equation*}

There is a natural decomposition of $\PS(G)$ into locally closed irreducible nonsingular
subvarieties, which we call \defterm{cellules}.  The decomposition is somewhat analogous to
the decomposition of a flag variety into Schubert cells.

\begin{defn}
Let $\A = \{A_1, \dots, A_s\}$ be a partition of $V$.  The \defterm{cellule} of $\A$ in 
$\PS(G)$ is the quasiprojective subvariety
	\begin{equation*}
	\PS_\A(G) = \{\Pic \in \PS(G) \st \Pic(v)=\Pic(w) ~\iff~ v \sim_\A w\}.
	\end{equation*}
\end{defn}

Unlike a Schubert cell, a cellule $\PS_\A(G)$ is not isomorphic to an affine space.  It is,
however, a smooth fiber bundle.  To see this, let $\Pic \in \PS_\A(G)$ and $e = vw \in E$.
If $v \sim_\A w$, then the set of lines in $\Pp^2$ through
$\Pic(v)=\Pic(w)$ is isomorphic to $\Pp^1$, and $\Pic(e)$ may take any value in that set.
If on the other hand $v \not\sim_\A w$, then $\Pic(e)$ is determined
uniquely by $\Pic(v)$ and $\Pic(w)$.  Therefore, putting
	\begin{equation*}
	q = \big\vert \{vw \in E \st v \sim_A w\} \big\vert
	\end{equation*}
and
	\begin{equation*}
	U = \big\{ (p_1, \dots, p_s) \in (\Pp^2)^s \st i \neq j ~\implies~ p_i \neq p_j 
	\big\},
	\end{equation*}
we see that $\PS_\A(G)$ has the bundle structure
	\begin{equation} \label{fiberbundle}
	\begin{array}{ccc}
	(\Pp^1)^q & \to & \PS_A(G) \\
	&& \downarrow \\
	&& U
	\end{array} \quad .
	\end{equation}
In particular, this yields the useful formula
	\begin{equation}\label{dimcellule}
	\dim \PS_\A(G) = 2s + \big\vert \{vw \in E \st v \sim_\A w\} \big\vert.
	\end{equation}

\begin{defn}
Let $G=(V,E)$.  A picture $\Pic \in \PS(G)$ is called \defterm{generic} if no two of the
points $\Pic(v)$ coincide.  The \defterm{discrete cellule} $\VO(G)$ is defined as the set of
all generic pictures of $G$.  Note that $\VO(G) = \PS_\D(G)$, where $\D=\D_V$ is the discrete
partition of $V$ (the partition into singleton sets).  The \defterm{picture variety} of $G$
is defined as the closure of the discrfete cellule:
	\begin{equation*}
	\PV(G) = \overline{\VO(G)}.
	\end{equation*}
This is an irreducible component of $\PS(G)$.  By (\ref{dimcellule}), we have
	\begin{equation*}
	\dim \PV(G) = \dim \VO(G) = 2|V(G)|.
	\end{equation*}
The \defterm{affine picture variety} of $G$ is defined as
	\begin{equation*}
	\PVA(G) = \PV(G) \cap \PSA(G).
	\end{equation*}
\end{defn}

\begin{remark} \label{blowup-rmk}
For $G=(V,E)$ and $W \subset V$, define the \defterm{coincidence locus} of $W$ as
	\begin{equation} \label{coinc-locus}
	C_W(G) = \{\Pic \in \PS(G) \st \Pic(v)=\Pic(w) ~\text{for all}~
	v,w \in W\}.
	\end{equation}
Let $G_0$ be the graph with vertices $V=V(G)$ and no edges. We may regard $\PV(G)$ as the
simultaneous blowup of $(\Pp^2)^n = \PS(G_0)$ along the coincidence loci $C_{\{e\}}(G)$ for
all edges $e$.  Indeed, the further blowup of $(\Pp^2)^n$ along all $C_W$, where $W \subset V$
is connected, is an instance of the ``wonderful model of subspace arrangements'' of DeConcini
and Procesi \cite{DP95}.  This blowup is a desingularization of $\PV(G)$. When $G$ is the
complete graph $K_n$, this is the ``compactification of configuration space'' of Fulton and
MacPherson \cite{FM94}.
\end{remark}

Note that the only cellule which is closed in $\PS(G)$ is the \defterm{indiscrete cellule}
$\PS_\I(G)$, where $\I=\I_V$ is the indiscrete partition of $V$ (the partition with just one
block).

\begin{example}
Let $G = K_2$.  Denote by $\D$ and $\I$ respectively the discrete and indiscrete partitions
of $V=V(G)=\{1,2\}$.  The picture space $\PS(K_2)$ is the blowup of $\Pp^2 \x \Pp^2$ along
the diagonal
	\begin{equation*}
	\{(p_1,p_2) \in \Pp^2 \x \Pp^2 \st p_1 = p_2\}.
	\end{equation*}
The blowup map
	\begin{equation*}
	\PS(K_2) \to \Pp^2 \x \Pp^2
	\end{equation*}
is just the projection on the vertex coordinates.  The exceptional divisor of the blowup
is the indiscrete cellule, which has dimension 3.  Since there are no 
partitions of $V$ other than $\D$ and $\I$, the complement of $\PS_\I(K_2)$ is $\VO(K_2)$, 
which has dimension 4 and is dense in $\PS(K_2)$.  Thus $\PV(K_2) = \PS(K_2)$.
\end{example}

\begin{example}
In general, the picture space $\PS(G)$ is not irreducible.  The first example, and in many
ways the fundamental one, is the graph $K_4$.  Denote by $\D$ and $\I$ respectively the
discrete and indiscrete partitions of $\{1,2,3,4\}$.  The data for a picture in $\PS_\I(K_4)$
consists of six lines, each passing through one free point in the plane, so
	\begin{equation*}
	\dim \PS_\I(K_4) = 8 = \dim \VO(K_4)
	\end{equation*}
(this follows also from (\ref{dimcellule})).  In particular, $\PS_\I(K_4)$ is too big to be
contained in the closure of $\VO(K_4)$.  Hence $\PV(K_4) \neq \PS(K_4)$.  In fact, the
irreducible components of $\PS(K_4)$ are precisely $\PV(K_4)$ and $\PS_\I(K_4)$ (see 
Lemma~\ref{circuit-component-lemma}).
\end{example}

We will soon see that the polynomials defining $\PVA(G)$ as a subvariety of $\PSA(G)$ involve
only the variables $m_e$.  In order to study these polynomials in isolation, we define a third
type of graph variety.  As before, we identify $\Aa^2$ with an open affine subset of $\Pp^2$.

\begin{defn} \label{defslopevar}
Let $U$ be the (dense) set of pictures $\Pic \in \PV(G)$ such that no $\Pic(e)$ is the line at
infinity.  Accordingly, for each $e$, $\Pic(e) \cap \Aa^2$ is an affine line of the form
	\begin{equation*}
	\big\{ (x,y) \st \alpha_e x + \beta_e y = 1 \big\},
	\end{equation*}
with a well-defined ``slope'' $[\alpha_e : \beta_e] \in \Pp^1$.  Forgetting all the 
data of $\Pic$ except the slopes gives a map
	\begin{equation} \label{natural-surjection}
	\phi: U \to (\Pp^1)^r.
	\end{equation}
We define the \defterm{slope variety} $\SV(G)$ as the image of $\phi$.  An element of $\SV(G)$
is called a \defterm{slope picture} of $G$.  If $a_1^e \neq 0$ for all $e$, then we have an
\defterm{affine slope picture}.  Setting $m_e = a_0^e/a_1^e$, we may regard an affine slope
picture as a point $\mm = (m_e \st e \in E)$ of $\Spec R_G$.  The algebraic set $\SVA(G)$ of
all affine slope pictures is called the \defterm{affine slope variety} of $G$.  $\SVA(G)$ may
also be defined as the projection of $\PVA(G)$ on the coordinates $m_e$; since $\PVA(G)$ is
irreducible, so is $\SVA(G)$.
\end{defn}

\begin{remark}
Restricting $\phi$ to $\PVA(G)$ produces a surjective map
	\begin{equation} \label{VAtoSA}
	\phi: \PVA(G) \to \SVA(G).
	\end{equation}
Note that every fiber of $\phi$ has dimension at least 3, because translation and
scaling do not affect slopes of lines.
\end{remark}

We will show that the same ideal of $R_G$ cuts out $\SVA(G)$ set-theoretically as a subvariety
of $\Spec R_G$, and $\PVA(G)$ as a subvariety of $\PSA(G)$.  To study this ideal, we use tools
from the theory of combinatorial rigidity.

\section{Combinatorial Rigidity Theory} \label{rigidity-section}

The behavior of graph varieties is governed in various ways by a certain combinatorial object,
the \defterm{generic rigidity matroid}.  Accordingly, we begin this section by sketching the
elements of rigidity theory, collecting several facts which we will need later.  (Our
treatment here is necessarily brief; for a detailed exposition, see \cite{GSS93} or
\cite{Whi96}.)  The main new result of this section, Theorem~\ref{matroid-thm}, describes the
fundamental connection between the purely combinatorial theory of rigidity and the geometry of
graph varieties.

Let $G=(V,E)$ be a connected graph, and $\Pic$ a generic picture of $G$. For the sake of easy
visualization, we may take $k=\Rr$ for the moment (the requirement that $k$ be algebraically
closed is not needed for the notion of rigidity).  Imagine a physical model of $\Pic$ in which
the vertices and edges are represented by ``ball joints'' and ``rods'' respectively.  The rods
are considered to be fixed in length, but are free to rotate about the joints in the plane of
the picture. Intuitively, $G$ is \defterm{length-rigid}, or simply \defterm{rigid}, if the
physical realization of any generic picture of $G$ ``holds its shape.'' More precisely, $G$ is
rigid if the distance between any two vertices in a generic picture is determined by the
lengths of the edges in $E$, up to finitely many possibilities.  (This property is called
``generic rigidity'' in \cite{GSS93}, as distinguished from other types of rigidity which we
will not need here.)

For instance, let $G$ be the $4$-gon on vertices $1,2,3,4$.  $G$ is not rigid, since there are
infinitely many incongruent rhombuses with equal side lengths, as shown in
Figure~\ref{quadrilateral-figure}.

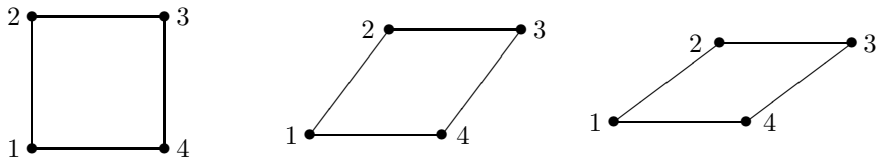
\begin{figure}
\begin{picture}(330,70)
\put( 3, 0){\makebox(0,0){1}}
\put( 3,50){\makebox(0,0){2}}
\put(67,50){\makebox(0,0){3}}
\put(67, 0){\makebox(0,0){4}}
\put(10, 0){\makebox(0,0){$\bullet$}}
\put(10,50){\makebox(0,0){$\bullet$}}
\put(60, 0){\makebox(0,0){$\bullet$}}
\put(60,50){\makebox(0,0){$\bullet$}}
\put(10, 0){\line(1,0){50}}
\put(10,50){\line(1,0){50}}
\put(10, 0){\line(0,1){50}}
\put(60, 0){\line(0,1){50}}
\put(108, 5){\makebox(0,0){1}}
\put(137,45){\makebox(0,0){2}}
\put(202,45){\makebox(0,0){3}}
\put(173, 5){\makebox(0,0){4}}
\put(115, 5){\makebox(0,0){$\bullet$}}
\put(165, 5){\makebox(0,0){$\bullet$}}
\put(145,45){\makebox(0,0){$\bullet$}}
\put(195,45){\makebox(0,0){$\bullet$}}
\put(115, 5){\line(3,4){30}}
\put(165, 5){\line(3,4){30}}
\put(115, 5){\line(1,0){50}}
\put(145,45){\line(1,0){50}}
\put(223,10){\makebox(0,0){1}}
\put(261,40){\makebox(0,0){2}}
\put(327,40){\makebox(0,0){3}}
\put(289,10){\makebox(0,0){4}}
\put(230,10){\makebox(0,0){$\bullet$}}
\put(280,10){\makebox(0,0){$\bullet$}}
\put(270,40){\makebox(0,0){$\bullet$}}
\put(320,40){\makebox(0,0){$\bullet$}}
\put(230,10){\line(4,3){40}}
\put(280,10){\line(4,3){40}}
\put(230,10){\line(1,0){50}}
\put(270,40){\line(1,0){50}}
\end{picture}
\caption{The $4$-gon is not rigid}
\label{quadrilateral-figure}
\end{figure}

However, the graph $G' = (V,E \cup \{24\})$ is rigid, because a generic affine picture of
$G'$ is determined by the lengths of its edges, up to isometries of $\Aa^2$ and finitely
many possibilities.

\begin{defn}
The \defterm{length-rigidity matroid} $\LRM{V}$ (called in \cite{GSS93}) the 
\defterm{2-dimensional generic rigidity matroid}) is the algebraic dependence 
matroid on the squares of lengths of edges
	\begin{equation*}
	(x_v-x_w)^2+(y_v-y_w)^2, \qquad v,w \in V.
	\end{equation*}
\end{defn}

We may regard $\LRM{V}$ as a matroid on $K(V)$, associating the above polynomial with the edge
$vw$. Accordingly, we say that a set of edges is independent in $\LRM{V}$, or
\defterm{rigidity-independent}, if and only if the corresponding set of squared lengths is
algebraically independent over ${\mathbb Q}$.  Thus an edge set $E$ is rigid if and only if
$E$ is a spanning set of $\LRM{V}$.

A fundamental result of rigidity theory is the following characterization of the independent
sets and bases of $\LRM{V}$, originally due to G.~Laman~\cite[Theorem 4.2.1]{GSS93}. An edge
set $E \subset K(V)$ is rigidity-independent if and only if it satisfies \defterm{Laman's
condition}:
	\begin{equation} \label{lamanone}
	|F| \leq 2 |V(F)|-3 ~\text{for all}~ F \subset E.
	\end{equation}
Furthermore, a rigidity-independent set $E$ is a basis of $\LRM{V}$ if and only if
	\begin{equation} \label{lamantwo}
	|E| = 2|V|-3.
	\end{equation}
In addition, $E$ is a \defterm{rigidity circuit}---a minimal dependent set of
$\LRM{V}$---if and only if $|E|=2|V(E)|-2$ and every proper subset $F$ of $E$ satisfies
Laman's condition (\ref{lamanone}) \cite[Theorem 4.3.1]{GSS93}.

The rigidity circuits (called ``rigidity cycles'' in \cite{GSS93}) may be described another
way.  Define a \defterm{rigidity pseudocircuit} to be an edge set $E$ equal to the
edge-disjoint union of two spanning trees of $V(E)$.  Then a rigidity circuit is a minimal
rigidity pseudocircuit \cite[Lemma 4.9.3 and Theorem 4.9.1]{GSS93}.

\begin{example} \label{rcircuitex}
Let $r \geq 3$.  The \defterm{$r$-wheel} is the graph $W_r$ with vertices 
	\begin{equation*}
	\{v_0, v_1, \dots, v_r\}
	\end{equation*}
and edges
	\begin{equation*}
	\{v_1v_2, v_2v_3, \dots, v_rv_1\} \cup \{v_0v_1, v_0v_2, \dots, v_0v_r\}.
	\end{equation*}
For all $r \geq 3$, $W_r$ is a rigidity circuit \cite[Exercise 4.13]{GSS93}. (In
fact, $W_3 \isom K_4$ and $W_4$ are the only rigidity circuits on 5 or fewer vertices.)  On
the other hand, let $G'$ be the graph given in Figure~\ref{pseudocircuit-figure}.
This is a rigidity pseudocircuit, since its edges are the disjoint union of the spanning trees
$\{12, 23, 34, 45\}$ and $\{13, 14, 24, 35\}$, but it is not a rigidity circuit since it
contains $K_4$ as a proper subgraph.
\end{example}

\begin{figure}
\begin{picture}(160,100)
\put(  0, 40){\makebox(0,0){1}}
\put( 70, 40){\makebox(0,0){2}}
\put(110,100){\makebox(0,0){3}}
\put(110,  0){\makebox(0,0){4}}
\put(160, 50){\makebox(0,0){5}}
\put(  0, 50){\makebox(0,0){$\bullet$}}  
\put( 75, 50){\makebox(0,0){$\bullet$}}  
\put(100,100){\makebox(0,0){$\bullet$}}  
\put(100,  0){\makebox(0,0){$\bullet$}}  
\put(150, 50){\makebox(0,0){$\bullet$}}  
\put( 75, 50){\line(-1, 0){ 75}}  
\put(100,100){\line(-2,-1){100}}  
\put(100,  0){\line(-2, 1){100}}  
\put(100,100){\line(-1,-2){ 25}}  
\put(100,  0){\line(-1, 2){ 25}}  
\put(100,  0){\line( 0, 1){100}}  
\put(100,100){\line( 1,-1){ 50}}  
\put(100,  0){\line( 1, 1){ 50}}  
\end{picture}
\caption{A pseudocircuit which is not a circuit}
\label{pseudocircuit-figure}
\end{figure}
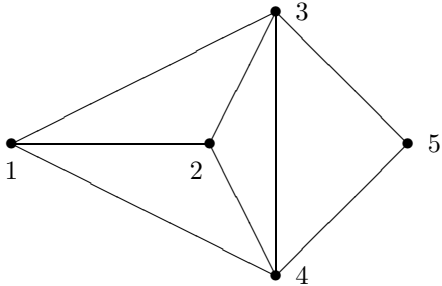

\begin{defn}
Let $G=(V,E)$ be a rigidity pseudocircuit.  A \defterm{coupled spanning tree} of $G$ is an
edge set $T \subset E$ such that both $T$ and $T'=E \sm T$ are spanning trees of $V$.  The
set of coupled spanning trees of $G$ is denoted by $\Coup(G)$.  The pair $T,T'$ is called a
\defterm{2-tree decomposition} of $E$ (or of $G$).
\end{defn}

The coupled spanning trees of a rigidity circuit play a fundamental role in describing
the equations which define $\PVA(G)$ and $\SVA(G)$.

Our local affine coordinates on $\PSA(G)$ measure the slopes of edges rather than their
lengths, leading to an alternate notion of rigidity.

\begin{defn}
The \defterm{slope-rigidity matroid} $\SRM{V}$ on $K(V)$ is the algebraic dependence matroid
on the rational functions
	\begin{equation*}
	m_{vw} = \frac{y_w-y_v}{x_w-x_v}.
	\end{equation*}
\end{defn}

\begin{theorem} \label{matroid-thm}
Let $G = (V,E)$.  The following are equivalent:

\noindent {\rm (i)} $E$ is independent in $\LRM{V}$;

\noindent {\rm (ii)} $E$ is independent in $\SRM{V}$;

\noindent {\rm (iii)} $\SVA(G) = \Spec R_G$ ($\isom \Aa^{|E|}$);

\noindent {\rm (iv)} $\PV(G)=\PS(G)$.
\end{theorem}

\begin{proof}
Let $n=|V|$ and $r=|E|$.

\noindent (i) $\implies$ (iv): Since $\PS(G)$ is defined locally by $2r$ equations among
$2n+2r$ coordinates, we have
	\begin{equation*}
	\dim X \geq 2n = \dim \PV(G)
	\end{equation*}
for every irreducible component $X$ of $\PS(G)$.  Therefore $\VO(G)$ is dense in 
$\PS(G)$ if and only if every nondiscrete cellule has dimension $< 2n$.

Suppose $E$ is rigidity-independent, hence satisfies Laman's condition (\ref{lamanone}).  
Let $\A$ be a partition of $V$ which is not the discrete partition. The blocks of $\A$ may
be numbered $A_1, \dots, A_s$ so that
	\begin{equation*}
	\begin{array}{ll}
	|A_i| = 1 & \text{for } 1 \leq i \leq t, \\
	|A_i| > 1 & \text{for } t+1 \leq i \leq s. \\
	\end{array}
	\end{equation*}
We may rewrite the cellule dimension formula (\ref{dimcellule}) as
	\begin{equation*}
	\dim \PS_\A(G) ~=~ 2s ~+~ \sum_{i=1}^s |K(A_i) \cap E|.
	\end{equation*}
If $i \leq t$, then $K(A_i) = \emptyset$, while if $i > t$, then $|K(A_i) \cap E|
\leq 2|A_i|-3$ by Laman's condition (\ref{lamanone}).  Hence
	\begin{equation} \label{dimcellbound}
	\begin{split}
	\dim \PS_\A(G) &\leq~ 2s + \sum_{i=t+1}^s (2|A_i|-3) ~=~ 2s + (2(n-t)-3(s-t)) \\
	&=~ 2n-s+t ~<~ 2n.
	\end{split}
	\end{equation}
as desired.

\noindent (iv) $\implies$ (iii): No nonzero element of $R_G$ vanishes on $\PSA(G)$, since the
projection of the indiscrete cellule $\PSA_\I(G)$ on the second factor in (\ref{grg}) is
surjective.  On the other hand, every element of $R_G$ that vanishes on $\SVA(G)$ vanishes on
$\PVA(G)$.  We conclude that if (iii) fails, then (iv) fails as well.

\noindent (iii) $\implies$ (ii): This is essentially the definition of the slope-rigidity matroid.

\noindent (ii) $\implies$ (i): Suppose that $E$ is independent in $\SRM{V}$.  Let $F \subset
E$, and let $H$ be the graph $(V(F),F)$.  Then $\dim \PVA(H) = 2|V(F)|$, and all fibers of the
canonical surjection $\PVA(H) \to \SVA(H)$ have dimension $\geq 3$ (since translation and
scaling do not affect slope), whence
	\begin{equation} \label{ineqone}
	\dim \SVA(H) \leq 2|V(F)|-3.
	\end{equation}
On the other hand, $F$ is independent in $\SRM{V}$, so $\{m_f \st f \in F\}$ is algebraically 
independent, and these variables form a system of parameters for $\SVA(H)$, so
	\begin{equation} \label{ineqtwo}
	|F| \leq \dim \SVA(H).
	\end{equation}
Together, (\ref{ineqone}) and (\ref{ineqtwo}) imply Laman's condition 
(\ref{lamanone}) for $E$.
\end{proof}

We have recovered the following fact.
\footnote{During the refereeing process, the author learned that his result had already 
appeared in the literature: see, e.g., \cite[Corollary 4.1.3]{Whi96}.  The author wishes 
to thank Walter Whiteley for helpful conversations on this topic.}

\begin{corollary}
For every vertex set $V$, the length-rigidity matroid $\LRM{V}$ and the slope-rigidity 
matroid $\SRM{V}$ are equal.
\end{corollary}

\section{The Equations Defining $\PVA(G)$} \label{main-section}

Let $G=(V,E)$ be a connected graph.  In this section, we explicitly construct an ideal $I=I_G$
defining the affine picture and affine slope varieties of $G$ set-theoretically.  The
generators of $I$ turn out to have an elegant combinatorial description: their terms enumerate
coupled spanning trees of the rigidity circuit subgraphs of $G$.

We begin with some computations which are most conveniently expressed in terms of the homology
of $G$, considered as a 1-dimensional simplicial complex.

A \defterm{directed edge} of $G$ is a symbol $\de{v}{w}$, where $vw \in E$.  An
\defterm{orientation} of an edge $e=vw$ is chosen by putting either $e = \de{v}{w}$ or $e =
\de{w}{v}$.  In addition, for $e = \de{v}{w} \in E$, we define
	\begin{equation*}
	x_f = x_w-x_v.
	\end{equation*}
In what follows, we fix an arbitrary orientation for each edge of $G$.
Let $C$ be the free $\Zz$-module on the directed edges of $G$, modulo the relations
	\begin{equation*}
	\de{w}{v} = -\de{v}{w}.
	\end{equation*}
Homologically, $C$ is the group of 1-chains.  The subgroup of 1-cycles is
	\begin{equation}
	Z = \left\{ \sum_e c_e e \in C \st \sum_e c_e x_e = 0 \right\}.
	\end{equation}
Note that $Z$ is generated by the cycles
	\begin{equation} \label{cycle}
	z(P) = \sum_{i=1}^s \de{v_i}{v_{i+1}}
	\end{equation}
where $P = (v_1, \dots, v_s, v_{s+1}=v_1)$ is a polygon of $G$.

The \defterm{support} of a chain is defined by
	\begin{equation}
	\supp \left( \sum_{e \in E} c_e e \right) = \left\{ e \in E \st c_e \neq 0 \right\}.
	\end{equation}
Note that if $\gamma \in Z$ and $\supp(\gamma) \subset T$ for some tree $T$, then $\gamma=0$.

Let $T$ be a spanning tree of $G$ and $S = E \sm T$.  Fix an orientation for each edge of $E$.
For each edge $e = \de{v}{w} \in S$, the edge set $T \cup \{e\}$ contains a unique polygon of
the form
	\begin{equation} \label{pte}
	P_T(e) = (v=v_1, \dots, v_s=w, v).
	\end{equation}
There is a corresponding cycle
	\begin{eqnarray}
	z_T(e) &=& \de{v_1}{v_2} + \dots + \de{v_{s-1}}{v_s} + \de{v_s}{v_1} \notag \\
	&=& -\de{v}{w} + \sum_{i=1}^{s-1} \de{v_i}{v_{i+1}} \notag \\
	&=& -e + \sum_{f \in T} c^T_{ef} ~ f \label{Tpolygon}
	\end{eqnarray}
where $c^T_{ef} \in \{0,1,-1\}$ for all $f$.  Note that for every spanning tree $T$ of $G$, 
the set
	\begin{equation*}
	\{ z_T(e) \st e \in T \}
	\end{equation*}
generates $Z$.  Indeed, if $\zeta = \sum_{e \in E} b_e e$ is a cycle, then
	\begin{equation*}
	\zeta' = \zeta + \sum_{e \in E \sm T} b_e z_T(e)
	\end{equation*}
is a cycle with support contained in $T$, so $\zeta' = 0$.

There is an injective map of $\Zz$-modules $C \to R'_G$ sending
	\begin{equation}
	\de{v}{w} \mapsto m_{vw}(x_v-x_w)
	\end{equation}
for all directed edges $\de{v}{w}$.  The image of $Z$ under this map contains
all polynomials $L(P)$ defined in (\ref{affinethree}). Therefore, for every
spanning tree $T$, the set $\left\{ L(P_T(e)) \st e \in S \right\}$
generates an ideal defining $\PSA(G)$ set-theoretically.

Let $e = \de{a}{b} \in S$, and let $P_T(e)$ be the polygon of (\ref{pte}).  Then
	\begin{equation} \label{lpte} 
	\begin{aligned}
	L(P_T(e)) ~&=~ \left( \sum_{i=1}^{s-1} m_{e_i} x_{e_i} \right) + m_e x_e \\
	&=~ -\sum_{f \in T} c^T_{ef} m_f x_f + m_e \sum_{f \in F} c^T_{ef} x_f \\
	&=~ \sum_{f \in T} c^T_{ef} (m_e-m_f) x_f~.
	\end{aligned}
	\end{equation}
Collecting the equations (\ref{lpte}) for all edges of $S$, we obtain a matrix equation
	\begin{equation} \label{main-matrix-eqn}
	\Big[ L(P_T(e)) \Big]_{e \in S} ~=~ M_T X_T ~=~ \left( D_S C_T - C_T D_T \right) X_T.
	\end{equation}
These matrices are defined as follows:
$M_T = [c^T_{ef}(m_e-m_f)]_{e \in S,\;f \in T}$;
$X_T$ is the column vector $\left[x_f\right]_{f \in T}$;
$C_T = [c^T_{ef}]_{e \in S,\;f \in T}$; and
$D_T$ (resp. $D_S$) is the diagonal matrix with entries $m_f$, $f \in T$ (resp. $m_e$, $e \in S$).
Accordingly, the equations (\ref{affinethree}) defining $\PS(G)$ are equivalent to
the single matrix equation
	\begin{equation} \label{mtxt}
	M_T X_T = 0.
	\end{equation}

\begin{example} \label{kfour}
Let $G = K_4$.  Orient each edge $ab \in E(G)$ as $\de{a}{b}$, where $a<b$. Let $T =
\{\de{1}{2}, \de{1}{3}, \de{1}{4}\}$.  For $e=\de{v}{w} \not\in T$, we have
	\begin{equation*}
	P_T(e) = \de{v}{w} + \de{w}{1} + \de{1}{v}
	= \de{v}{w} - \left( \de{1}{w} - \de{1}{v} \right),
	\end{equation*}
so the matrix $C_T$ is
	\begin{equation*}
	\begin{array}{crcl}
	&& f \\
	& \de{1}{2} & \de{1}{3} & \de{1}{4} \\ \\
	\begin{array}{cc} & \de{2}{3} \\ e & \de{2}{4} \\ & \de{3}{4} \end{array} &
	\left[ \begin{array}{c} -1 \\ -1 \\ 0 \end{array} \right. &
	\begin{array}{c} 1 \\ 0 \\ -1 \end{array} &
	\left. \begin{array}{c} 0 \\ 1 \\ 1 \end{array} \right]~,
	\end{array}
	\end{equation*}
the polynomials of (\ref{lpte}) are
	\begin{equation*}
	L\left(P_T\left(\de{a}{b}\right)\right) =
	m_{ab}(x_b-x_a) + m_{1b}(x_b-x_a) + m_{1b}(x_a-x_b)
	\end{equation*}
for $\de{a}{b} \in S = \{\de{2}{3},\de{2}{4},\de{3}{4}\}$, and the matrix of
(\ref{main-matrix-eqn}) is
	\begin{equation*} M_T X_T ~=~
	\left[
		\begin{array}{ccc}
		m_{12}-m_{23} & m_{23}-m_{13} & 0 \\
		m_{12}-m_{24} & 0 & m_{24}-m_{14} \\
		0 & m_{13}-m_{34} & m_{34}-m_{14}
		\end{array}
	\right]
	\left[
		\begin{array}{c}
		x_2-x_1 \\ x_3-x_1 \\ x_4-x_1
		\end{array}
	\right]~.
	\end{equation*}
\end{example}

\begin{lemma} \label{inverse-tree-lemma}
Let $G=(V,E)$ be a rigidity circuit, $T$ a coupled spanning tree, and $S = E \sm T$.  (Recall 
that $S$ is also a spanning tree of $G$, and that $M_T$, $C_T$, and $C_S$ are $(|V|-1) \x 
(|V|-1)$ square matrices.)  Then $C_T = C_S^{-1}$.
\end{lemma}

\begin{proof}
Replacing each edge $f$ on the right side of (\ref{Tpolygon}) with $f+z_S(f)$, we see that
	\begin{equation*}
	e - \sum_{f \in T} c^T_{ef} \sum_{g \in S} c^S_{fg} g ~=~
	e - \sum_{g \in S} g \sum_{f \in T} c^T_{ef} c^S_{fg} ~\in~ Z.
	\end{equation*}
This cycle is zero, because its support is contained in $S$.  Hence, for all $e,g \in S$,
	\begin{equation*}
	\sum_{f \in T} c^T_{ef} ~ c^S_{fg} ~=~ \delta_{eg},
	\end{equation*}
the Kronecker delta.  This is precisely the statement that $C_T = C_S^{-1}$.
\end{proof}

Somewhat more generally, if $G$ is a rigidity pseudocircuit and $T,U$ are two spanning trees
of $G$, then every polynomial in the set
	\begin{equation*}
	\{L(P_U(e)) \st e \in E \sm U\}
	\end{equation*}
may be expressed as an integer linear combination of the polynomials
	\begin{equation*}
	\{L(P_T(e)) \st e \in E \sm T\},
	\end{equation*}
and vice versa.  In addition, each of these sets is linearly independent, since, e.g., for $e 
\not\in T$, each variable $m_e$ appears in exactly one $L(P_T(e))$.  Therefore
	\begin{equation*}
	M_T X_T = B M_U X_U
	\end{equation*}
for some invertible integer matrix $B$.  In particular $\det B = \pm 1$, so the polynomial
$\det M_T$ is independent, up to sign, of the choice of $T$.  This motivates the following
definition:

\begin{defn}
Let $G=(V,E)$ be a connected graph with $|E|=2|V|-2$ (so that $M_T$ is a $(|V|-1) \x (|V|-1)$
square matrix) and let $T$ be a spanning tree.  The \defterm{tree polynomial} of $G$ is
defined up to sign as
	\begin{equation*}
	\TP(G) = \det M_T
	\end{equation*}
where $M_T$ is the matrix of (\ref{main-matrix-eqn}).
\end{defn}

It is immediate from the construction of $M_T$ that $\TP(G)$ is homogeneous of degree $|V|-1$.
The name ``tree polynomial'' is justified by the following theorem.  One more piece of
notation: to each edge set $F \subset E$, we associate the squarefree monomial
	\begin{equation*}
	m_F = \prod_{f \in F} m_f.
	\end{equation*}

\begin{theorem} \label{tree-poly-thm}
Let $G = (V,E)$ be a graph with $|V|=n$ and $|E|=2n-2$.
The tree polynomial of $G$ can be expressed as the sum
	\begin{equation*}
	\TP(G) = \sum_{T \in \Coup(G)} \sgn(T) m_T
	\end{equation*}
of signed squarefree monomials for coupled spanning trees
$T$, where $\sgn(T) \in \{1,-1\}$, and $\sgn(E \sm T)
= (-1)^{n-1} \sgn(T)$.  The tree polynomial vanishes
on both $\PVA(G)$ and $\SVA(G)$.
Moreover, $\TP(G)$ is nonzero if and only if G is a rigidity
pseudocircuit, and irreducible if and only if $G$ is a
rigidity circuit.  In particular, if $G'$ is a rigidity pseudocircuit 
subgraph of $G$, then $\TP(G')$ divides $\TP(G)$.
\end{theorem}

\begin{proof}
Fix a spanning tree $T$, so that $\TP(G) = \pm \det M_T$.  Let $S = E \sm T$.
For each edge $e$, if $e \in T$
then the variable $m_e$ appears in only one column of $M_T$, while if $e \in S$
then $m_e$ appears in only one row of $M_T$.  It follows that $\TP(G)$ is squarefree.

Each nonzero term in the determinant expansion of $\TP(G)$ is of the form
	\begin{equation*}
	(m_{1,1}-m_{1,2}) \dots (m_{n-1,1}-m_{n-1,2})
	\end{equation*}
where the $m_j$ are all distinct.  This may be expressed as a sum of binomials of the form
	\begin{equation*}
	m_F + (-1)^{n-1} m_{E \sm F},
	\end{equation*}
where each $F$ is a subset of $E$ of cardinality $n-1$.  It follows that $\sgn(F) = (-1)^{n-1}
\sgn(F)$ for all $F$.  In particular, $\sgn(S) = \det C_T$ by the definition of $M_T$
(\ref{main-matrix-eqn}).  If $S$ is also a tree, then by Lemma~\ref{inverse-tree-lemma}, $C_T$
is an invertible integer matrix, so its determinant is $\pm 1$.

Now suppose that $F \subset E$ has cardinality $n-1$, but is not a tree. We will show that
$\sgn(F)=0$. Let $A \subset F$ be a minimal set of edges such that $F \sm A$ is a forest; in
particular $A$ is nonempty. Let $T$ be a spanning tree of $G$ containing $F \sm A$; then
	\begin{equation*}
	T \cap (E \sm F) = T \sm (F \sm A) \neq \emptyset.
	\end{equation*}
Let $S = E \sm T$ ($\supset A$).  The matrix $C_T$ constructed in (\ref{main-matrix-eqn}) has
the property that $c^T_{ab}=0$ whenever $a \in A$ and $b \in T \cap (E \sm F)$, because the
unique circuit of $T \cup \{a\}$ is contained in $F \sm A \cup \{a\}$. Accordingly, for each
$a \in A$, every entry of the corresponding row of $M_T$ is either zero or of the form $\pm
(m_a-m_f)$, where $f \in F \sm A$.  In particular, no variable dividing $m_{E \sm F}$ appears
in that row.  Hence $\sgn(E \sm F)=0$, and $\sgn(F)=(-1)^{n-1}\sgn(E \sm F) = 0$ as well.  We
have obtained the desired equation
	\begin{equation*}
	\TP(G) ~= \sum_{T \in \Coup(G)} \sgn(T) m_T.
	\end{equation*}
By definition, the right side is nonzero if and only if $G$ is a rigidity pseudocircuit.

We next show that the tree polynomial vanishes on the affine picture and slope varieties.
Since the generic affine pictures are dense in $\PVA(G)$, and their image under the natural
surjection $\phi$~\ref{natural-surjection} is dense in $\SVA(G)$, it suffices to show that
$\TP(G)$ vanishes at each $\Pic \in \VO(G) \cap \PVA(G)$.  Indeed, $M_T(\Pic) X_T(\Pic) = 0$
and $X_T(\Pic) \neq 0$, so $\TP(G) = \det M_T$ vanishes at $\Pic$.

Suppose that $G$ contains a rigidity circuit $G'=(V',E')$ as a proper subgraph.  Let
$T'$ be a spanning tree of $G'$ and $T \supset T'$ a spanning tree of $G$.  Put $S = E \sm
T$ and $S' = E' \sm T'$.  Then the matrix $M_T$ has the form
	\begin{equation*}
	\left[ \begin{array}{cc}
	M_{T'} & 0 \\
	* & *
	\end{array} \right]
	\end{equation*}
where the $|V'|-1$ uppermost rows correspond to edges in $S'$ and the $|V'|-1$ 
leftmost columns correspond to edges in $T'$.  It follows that $\TP(G')$ is a proper divisor 
of $\TP(G)$.  In particular, if $\TP(G)$ is irreducible then $G$ is a rigidity circuit.

On the other hand, suppose that $G$ is a rigidity circuit and $\TP(G) = f_1 \cdot f_2$.
For every $e \in E$, we have
	\begin{equation*}
	\deg_{m_e}(\TP(G)) = \deg_{m_e}(f_1) + \deg_{m_e}(f_2) = 1,
	\end{equation*}
so $E$ may be expressed as a disjoint union $E_1 \cup E_2$, where $E_i = \{e \in E \st
\deg_{m_e}(f_i) = 1\}$. Let $G_i = (V,E_i)$.  Since $\SVA(G)$ is an irreducible variety,
either $f_1$ or $f_2$ must vanish on $\SVA(G)$; assume without loss of generality
that $f_1$ does so.  Then $f_1$ vanishes on $\SVA(G_1)$ as well via the natural surjection
$\SVA(G) \to \SVA(G_i)$.  By Theorem~\ref{matroid-thm}, $E_1$ must be rigidity-dependent.
But $E$ contains no proper rigidity-dependent subset, so we must have $E_1=E$.  Therefore
$E_2 = \emptyset$ and the factorization of $\TP(G)$ is trivial.
\end{proof}

\begin{example}
Let $G=K_4$.  Let $T$, $M_T$, $X_T$ be as in Example~\ref{kfour}.  There are two kinds of
spanning trees of $G$: paths $(a,b,c,d)$, and ``stars,'' such as $T$.  The paths are
coupled; the stars are not.  There are $4!/2 = 12$ paths, and the sign of a path is given by
the sign of the corresponding permutation in the symmetric group $S_4$, that is,
	\begin{equation*}
	\TP(K_4) = \det M_T = - \frac{1}{2} \sum_{\sigma \in S_4} sgn(\sigma)
	m_{\sigma_1 \sigma_2} m_{\sigma_2 \sigma_3} m_{\sigma_3 \sigma_4}~.
	\end{equation*}
On the other hand, if $G'$ is the graph of Example~\ref{rcircuitex} (a rigidity
pseudocircuit which is not a circuit), then
	\begin{equation*}
	\TP(G') = \pm \left( m_{35}-m_{45} \right) \TP(K_4).
	\end{equation*}
\end{example}

\begin{theorem} \label{defining-eqns-thm}
Let $G=(V,E)$ be a graph.  Let $I = I_G$ be the ideal of $R_G$ generated by all tree
polynomials $\TP(C)$, where $C$ is a rigidity circuit subgraph of $G$.  Then:
\begin{enumerate}
\item $\PVA(G)$ is the vanishing locus of $IR'_G$ in $\PSA(G)$.

\item $\SVA(G)$ is the vanishing locus of $I$ in $\Spec R_G$.
\end{enumerate}
\end{theorem}

\begin{proof}
We may assume without loss of generality that $G$ is connected, since every rigidity
circuit is connected and $\PVA(G)$ is the product of the picture varieties of its connected
components.  Let $n=|V|$, $r=|E|$.

Let $Y$ be the vanishing locus of $IR'_G$ in $\PSA(G)$.  For each rigidity circuit subgraph
$C$ of $G$, the tree polynomial $\TP(C)$ vanishes on $\PVA(C)$ by Theorem~\ref{tree-poly-thm},
so it vanishes on $\PVA(G)$ as well. Hence $\PVA(G) \subset Y$.

We now establish the reverse inclusion, proceeding by induction on $n$.  By
Theorem~\ref{matroid-thm}, there is nothing to prove when $E$ is rigidity-independent, in
particular when $n \leq 3$.

Let $\Pic \in Y \cap \PS_\A(G)$, where $\A = \{A_1, \dots, A_s\}$ is a partition of $V$ with
$s$ parts.  We wish to show that $\Pic \in \PVA(G)$.

\Case{1}{$s=n$}.  Here $\A$ is the discrete partition, so $\Pic \in \PVA(G)$ by definition.

\Case{2}{$2 \leq s \leq n-1$}.  For $1 \leq i \leq s$, define a subgraph
	\begin{equation*}
	G_i = (A_i,E \cap K(A_i))
	\end{equation*}
and let
	\begin{equation*}
	U = \bigcup_{\B \refineq \A} \PSA_\B(G) = \left\{ \Pic' \in \PSA(G) \st \Pic'(v) 
	\neq \Pic'(w) ~\text{if}~ v \not\sim_\A w \right\},
	\end{equation*}
an open subset of $\PSA(G)$ containing $\Pic$.  There is a natural open embedding
	\begin{equation}
	\theta: U \to \prod_{i=1}^s \PSA(G_i).
	\end{equation}
Note that $I_{G_i} \subset IR_G$ for all $i$.  By induction, $\PVA(G_i)$ is the
vanishing locus of $I_{G_i}$ in $\PSA(G_i)$.  Therefore
	\begin{equation*}
	\Pic \in \theta^{-1} \left( \prod_{i=1}^s \PVA(G_i) \right)~.
	\end{equation*}
This set is irreducible and contains $\VO(G)$ as 
an open, hence dense, subset.  Therefore $\Pic \in \PVA(G)$ as desired.

\Case{3}{$s=1$}.  That is, $\A$ is the indiscrete partition of $V$.
Fix a spanning tree $T$ of $G$ and let $M_T$ be the matrix defined in (\ref{main-matrix-eqn}).  
Recall that $M_T$ is an $(r-n+1) \x (n-1)$.  The rows and columns of $M_T$ are indexed 
by the edges of $E \sm T$ and $T$, respectively.  In addition, $\PSA(G)$ is defined by the
matrix equation (\ref{mtxt}), and $X_T(\Pic)$ is the zero matrix.

We claim that
	\begin{equation} \label{rankmt}
	\mathop{\rm rank}\nolimits M_T < n-1.
	\end{equation}
If $M_T$ has fewer than $n-1$ rows then there is nothing to prove.  Otherwise, let $M'$ be any
$(n-1) \x (n-1)$ submatrix $M'$ of $M_T$, with rows indexed by the elements of some edge set
$S \subset E \sm T$.  Then $|T \cup S|=2n-2$, so $|T \cup S|$ does not satisfy Laman's
condition (\ref{lamanone}) and must contain some rigidity circuit $C$.  By
Theorem~\ref{tree-poly-thm}, $\TP(C)$ divides $\det M'$, establishing (\ref{rankmt}).

It follows from (\ref{rankmt}) that the nullspace of $M_T(\Pic)$ contains a nonzero vector
$X'$.  For every $\lambda \in \fld$, we have $(M_T)(\lambda X') = 0$, so there is a picture
$\Pic_\lambda$ with the same slope coordinates as $\Pic$ and $x$-coordinates of vertices given
by $\lambda X'$.  The $\Pic_\lambda$ form an affine line in $Y$ with $\Pic_0=\Pic$.  
Moreover, if $\lambda \neq 0$, then $\Pic_\lambda \not\in \PS_\A(G)$, hence $\Pic_\lambda \in
\PVA(G)$ by the previous two cases.  Therefore $\Pic_0=\Pic \in \PVA(G)$ as well.

We now turn to the second assertion of the theorem.  Let $Z$ be the vanishing locus of $I$ in
$\Spec R_G$.  It is immediate from Definition ~\ref{defslopevar} that $Z \supset \SVA(G)$.  
Now suppose that $\mm \in Z$, i.e., $\mm$ is an affine slope picture at which all tree
polynomials vanish.  Fix a spanning tree $T$ of $G$ and let $X$ be a nullvector of the matrix
$M_T(\mm)$.  Together, $\mm$ and $X$ define an affine line in $\PSA(G)$; by part (i) of the
theorem, the line is contained in $\PVA(G)$.  Therefore $\mm \in \SVA(G)$.
\end{proof}

We have proven that
	\begin{equation}
	\PVA(G) \isom \Aa^1 \x \Spec R'_G/\sqrt{J_G}
	\end{equation}
and
	\begin{equation}
	\SVA(G) \isom \Spec R_G/\sqrt{I_G}
	\end{equation}
as reduced schemes, where $J_G = I_GR'_G+(L(P))$.  However, we do not yet know whether the
ideals $J_G$ and $I_G$ are radical.  In the special case that $G$ is a rigidity cycle, the
ideal $I_G$ is radical because it is principal, generated by the irreducible polynomial
$\TP(G)$.  We prove in a separate paper~\cite{JLM2} that $I_G$ is radical when $G$ is the
complete graph $K_n$.

\section{Further Geometric Properties of $\PS(G)$ and $\PV(G)$} \label{geom-section}

In this section, we use the algebraic results of the previous sections to prove certain
geometric facts about the picture space.  First, we give a combinatorial condition which
describes when one cellule of $\PS(G)$ is contained in the closure of another cellule.  
Using this result, we can give a complete combinatorial description of the irreducible
components of the picture space.  Second, we present an inductive criterion on $G$
which implies that $\PV(G)$ is Cohen-Macaulay; one consequence of this result is that
$\PV(G)$ is Cohen-Macaulay whenever $G$ is rigidity-independent.

\begin{defn}
Let $G = (V,E)$ be a graph, $F \subset E$, and $\A$ a partition of $V$.  We say that $\A$
\defterm{collapses} $F$ if all vertices of $V(F)$ are contained in the same block of $\A$.
\end{defn}

In this case, the equations defining $\PS(G)$ impose no restrictions on the slopes of the
lines $\Pic(e)$ for pictures $\Pic \in \PS_\A(G)$ and edges $e \in F$.

\begin{lemma} \label{circuit-component-lemma}
Let $G = (V,E)$ be a rigidity circuit.  Then the picture space $\PS(G)$ has two irreducible
components, both of dimension $2|V|$: the picture variety $\PV(G)$ and the indiscrete cellule
$\PS_\I(G)$.
\end{lemma}

\begin{proof}
The cellule dimension formula (\ref{dimcellule}) gives $\dim \PS_\I(G) = \dim \VO(G) = 2n$.
The indiscrete cellule is itself closed, so it is an irreducible component of $\PS(G)$. On the
other hand, if $\A$ is neither the discrete nor indiscrete partition, then $\dim \PS_\A(G) <
2n$ by (\ref{dimcellbound}) (since Laman's condition (\ref{lamanone}) holds for every proper
subset of a rigidity circuit).  Since all components of $\PS(G)$ have dimension at least $2n$,
we must have $\PS_\A(G) \subset \PV(G)$ for every such $\A$, which implies the desired result.
\end{proof}

\begin{theorem} \label{cellule-containment-thm}
Let $G = (V,E)$ be a graph, $\A$ a partition of $V$, and $\PS_\A(G)$ the corresponding cellule.
Then $\PS_\A(G) \subset \PV(G)$ if and only if no rigidity circuit of $G$ is collapsed by $\A$.
\end{theorem}

\begin{proof}
If $G$ is rigidity-independent, then this is immediate from Theorem~\ref{matroid-thm}, while
if $G$ is a rigidity circuit then the desired statement follows from
Lemma~\ref{circuit-component-lemma}. In general, by Theorem~\ref{defining-eqns-thm}, it is
enough to prove that for every rigidity circuit $C$ of $G$, $\TP(C)$ vanishes on $\PSA_\A(G)$
if and only if $\A$ does not collapse $C$.  One direction is immediate: if $\A$ collapses $C$,
then $\TP(C)$ does not vanish on $\PSA_\A(G)$ and consequently $\PSA_\A(G) \not\subset
\PVA(G)$.  On the other hand, suppose that $\A$ does not collapse $C$.  Consider the natural
map $\PSA_\A(G) \to \PSA(C)$.  The image of this map does not intersect the indiscrete cellule
of $\PSA(C)$.  By Lemma~\ref{circuit-component-lemma}, $\TP(C)$ vanishes on the image, hence 
on $\PS_\A(G)$.
\end{proof}

Given a graph $G=(V,E)$ and a partition $\A = \{A_1, \dots, A_s\}$ of $V$, let $G/\A$ denote
the graph whose vertices are the blocks of $\A$ and whose edges are
	\begin{equation*}
	\big\{ \ue{A_i}{A_j} \st vw \in E ~\text{for some}~ v \in A_i, ~w \in A_j \big\}.
	\end{equation*}
Also, if $\A$ and $\B$ are partitions of $V$ with $\A \refineq \B$, then we write
$\B/\A$ for the partition on the blocks of $\A$ setting two blocks equivalent if both are 
subsets of the same block of $\B$.

\begin{theorem} \label{cellule-closure-thm}
Let $G = (V,E)$ be a graph, $\A$ and $\B$ partitions of $V$, and $\PS_\A(G)$ and $\PS_\B(G)$ 
the corresponding cellules.  Then $\PS_\B(G) \subset \overline{\PS_\A(G)}$ if and only if the 
following conditions hold:

(a) $\A \refineq \B$;

(b) No rigidity circuit of $G/\A$ is collapsed by $\B/\A$; and

(c) If $A_i$ and $A_j$ are distinct blocks of $\A$ contained in the same block of $\B$, then
$E$ contains at most one edge between $A_i$ and $A_j$ (i.e., with one endpoint in each set).

Consequently, the irreducible components of $\PS(G)$ are exactly the subvarieties
$\overline{\PS_\A(G)}$, where $\A$ is maximal with respect to the partial order just described.
\end{theorem}

\begin{proof}
Suppose that $\PS_\B(G) \subset \overline{\PS_\A(G)}$.  If $v \sim_\A w$, the equation
$\Pic(v)=\Pic(w)$ holds on $\PS_\B(G)$, so $v \sim_\B w$, establishing (a). For each
rigidity circuit $C$ of $G/\A$, the function $\TP(C)$ vanishes on $\PS_\B(G)$, so $\B/\A$
cannot collapse $C$.  Finally, if $A_i$ and $A_j$ are contained in the same block of $\B$
and $E$ contains two distinct edges $e,e'$ between $A_i$ and $A_j$, then the equation
$\Pic(e)=\Pic(e')$ holds on $\PS_\A(G)$ but not on $\PS_\B(G)$, a contradiction.

Now suppose that conditions (a), (b) and (c) hold. It is harmless to replace $\PS_\A(G)$ and
$\PS_\B(G)$ with the affine cellules
	\begin{align*}
	\PSA_\A(G) ~&=~ \PS_\A(G) \, \cap \, \PSA(G), \\
	\PSA_\B(G) ~&=~ \PS_\B(G) \, \cap \, \PSA(G).
	\end{align*}
Let $E' = \{vw \in E \st v \sim_\A w\}$ (that is, the edges
corresponding to ``free lines'' in $\PS_\A(G)$).  Let $U = \Aa^{|E'|}$ and
	\begin{equation*}
	Z = \bigcup_{\B \refinedbyeq \A} \PSA_\B(G) =
	\left\{\Pic \in \PSA(G) \st \Pic(v)=\Pic(w) ~\text{if}~ v \sim_\A w\right\}.
	\end{equation*}
Observe that the data for an affine picture $\Pic \in Z$ is the same as that
describing a picture of $G/\A$ together with the slopes of the lines $\Pic(e)$ for
$e \in E'$.  Hence we have an isomorphism 
	\begin{equation*}
	\psi:~ Z \;\xrightarrow{\isom}\; \PSA(G/\A) \x U,
	\end{equation*}
Restricting $\psi$ to the cellules under consideration, we have a commutative diagram of
quasiaffine varieties:
	\begin{equation*}
	\begin{array}{ccccccc}
	\PSA_\A(G) & \subset & \overline{\PSA_\A(G)} & \subset & Z & \supset & \PSA_\B(G) \\
	\downarrow  && \downarrow && \downarrow && \downarrow \\
	\VAO(G/\A) \x U & \subset & \PVA(G/\A) \x U & \subset & \PSA(G/\A) \x U
		& \supset & \PSA_{\B/\A}(G/\A) \x U 
	\end{array}
	\end{equation*}
where the vertical arrows are isomorphisms.  This implies that
	\begin{equation*}
	\PSA_\B(G) \subset \overline{\PSA_\A(G)} ~\iff~ \PSA_{\B/\A}(G/\A) \subset \PVA(G/\A)
	\end{equation*}
which is in turn equivalent to condition (b) by Theorem~\ref{cellule-containment-thm}.
\end{proof}

\begin{remark}
The notion of a pseudocircuit may be extended to multigraphs: a multigraph $(V,E)$ is called
a \defterm{pseudocircuit} if $|E|=2|V|-2$ and $|F| \leq 2|V(F)|-2$ for all $\emptyset
\neq F \subset E$ \cite[p.~118]{GSS93}.  For instance, a double edge is a pseudocircuit.  
In the previous theorem, we may consider $G/\A$ as a multigraph, in which the multiplicity
of an edge $\ue{A_i}{A_j}$ is the number of edges in $E$ with one endpoint in each of $A_i$
and $A_j$. Then conditions (b) and (c) together are equivalent to the single condition that
$\B/\A$ collapse no multigraph pseudocircuit of $G/\A$.
\end{remark}

We next consider the Cohen-Macaulay property.  Our main tool is the fact that if $X$ is a
Cohen-Macaulay scheme and $Z$ is a ``strongly Cohen-Macaulay'' subscheme of $X$, then the
blowup of $X$ along $Z$ is Cohen-Macaulay \cite[Theorem 4.2]{Hu83} (see also \cite{SV80}).  
In particular, a local complete intersection subscheme of a Cohen-Macaulay scheme is
strongly Cohen-Macaulay.

\begin{lemma} \label{blowup-lemma}
Let $G = (V,E)$, $e = vw \in E$, and $H = (V, E \sm \{e\})$.  Suppose that $\PV(H)$ is
Cohen-Macaulay and that $\PV(H) \cap \PS_\A(H)$ has codimension $\geq 2$ in $\PV(H)$ for all
partitions $\A$ of $V$ with $v \sim_\A w$.  Then $\PV(G)$ is Cohen-Macaulay.
\end{lemma}

\begin{proof}
Let $Z$ be the (possibly non-reduced) intersection $\PV(H) \cap C_e(H)$, where
$C_e(H)=C_{\{e\}}(H)$ is the coincidence locus defined in (\ref{coinc-locus}).  $Z$ is defined
in local affine coordinates by two equations, namely $x_v=x_w$ and $y_v=y_w$, so each of its
components has codimension $\leq 2$.  On the other hand, $C_{\{e\}}(H)$ is set-theoretically
the union of cellules $\PS_A(H)$ with $v \sim_\A w$.  Therefore
	\begin{equation*}
	\codim Z \geq \codim C_{\{e\}}(H) \geq 2.
	\end{equation*}
In particular, $Z$
is a local complete intersection in $\PV(H)$, and $\PV(G)$ is the blowup of $\PV(H)$ along $Z$,
so $\PV(G)$ is Cohen-Macaulay.
\end{proof}

\begin{proposition} \label{CM-prop}
Let $G = (V,E)$, $e = vw \in E$, and $H = (V,E \sm \{e\})$.  If $\PV(H)$ is
Cohen-Macaulay and $e$ is not contained in any rigidity circuit subgraph of $G$, then
$\PV(G)$ is Cohen-Macaulay.
\end{proposition}

\begin{proof}
Let $\A$ be a partition of $V$ with $v \sim_\A w$.  The cellule $\PV_\A(G) = \PS_\A(G)  
\cap \PV(G)$ has codimension $\geq 1$ in $\PV(G)$.  Since no rigidity circuit contains $e$,
the equations defining $\PV_\A(G)$ impose no constraints on the line $\Pic(e)$.  Therefore
	\begin{equation*}
	\PV_\A(G) \isom \PV_\A(H) \x \Pp^1.
	\end{equation*}
In particular $\PV_\A(H)$ has codimension $\geq 2$ in $\PV(H)$, since $\dim \PV(G) =
\dim \PV(H) = 2|V|$.  Thus $\PV(G)$ is Cohen-Macaulay by Lemma~\ref{blowup-lemma}.
\end{proof}

\begin{theorem} \label{CM-thm}
Let $G=(V,E)$.  If $G$ is rigidity-independent, then $\PV(G)$ is Cohen-Macaulay.
\end{theorem}

\begin{proof}
If $E = \emptyset$, the result is trivial since $\PV(G) \isom (\Pp^2)^{|V|}$.  
Otherwise, we add one edge at a time, applying Proposition~\ref{CM-prop} at each stage.
\end{proof}

\end{document}